\documentclass[12pt]{article}
\usepackage{amsmath,amstext, verbatim,remreset}
\usepackage[T2A]{fontenc}
\usepackage[cp1251]{inputenc}
\usepackage[english]{babel}
\usepackage{amsfonts,amssymb, amsmath}

\usepackage{mathrsfs}
\renewcommand{\leq}{\leqslant}
\renewcommand{\geq}{\geqslant}

\makeatletter
\renewcommand{\@begintheorem}[2]{\vskip 0.1 cm \noindent {\bf #1\ #2.\ }\begingroup\it}
\renewcommand{\@endtheorem}{\endgroup \vskip 0.1 cm}
\makeatother

\begin{document}

\begin{center}
{\bf \Large Kolmogorov problem on the class of multiply monotone functions}

\bigskip

O. V. Kovalenko

\bigskip

Oles Gonchar Dnipropetrovsk National University
\end{center}
{\small \it \noindent Necessary and sufficient conditions for positive numbers
$M_{k_1}, M_{k_2},  M_{k_3}, M_{k_4}$, $0 = k_1 < k_2<k_2\leq r-2$, $k_4=r$, to guarantee the existence of an $r-1$-monotone function defined
on the negative half-line and such that $\|x^{(k_i)}\| = M_{k_i}$, $i=1,2,3,4$ were found. \hfill} \vskip 1.5mm

\bigskip

{\bf 1. Notations. Statement of the problem. Known results.} 
Denote by $G$ real line $\mathbb R=(-\infty , \infty )$ or non-positive semi-line $\mathbb R_-=(-\infty , 0]$. Let  $L_{\infty}(G)$ denote the space of all measurable essentially bounded functions  $x:G\to \mathbb R$ with usual norm $\| \cdot \|=\|
\cdot \|_{L_{\infty}(G)}$. For $r\in \mathbb N$ denote by $L^r_{\infty} (G)$  the space of functions $x:G\to \mathbb R$, that have locally absolute continuous derivative of order $r-1$, $x^{(0)}=x$, and such that $x ^ {(r)} \in L_{\infty} (G)$. Let $L^r_{\infty,\infty}(G)=L^r_{\infty}(G)\cap L_{\infty}(G)$. 

For real number $t\in \mathbb R$ let $t_+ := \max\{t,0\}$.

A. N. Kolomogorov (see~{[\ref{Kolmogorov}]}) stated the following problem:

{\bf Kolmogorov problem}

{\it Let some class of functions $X\subset L^r_{\infty,\infty}(G)$ and arbitrary system of  $d$ positive integers 
$0 = k_1 < k_2 < ... < k_{d}=r$ are given.
The problem is to find necessary and sufficient conditions on the system of positive numbers
$$M_{k_1}, M_{k_2}, ..., M_{k_{d}}$$
to guarantee the existence of a function $x \in X$ such that $$\|x^{(k_i)}\|=M_{k_i},\; i=1,...,d.$$
}

In~{[\ref{Kolmogorov}] A. N. Kolomogorov solved this problem in the case $d=3$, $X= L^r_{\infty,\infty}(\mathbb R)$ (partial cases follow from works of Hadamard~[\ref{had}] and Shilov~[\ref{shilov}]). 
He showed that for three positive numbers $M_0, M_k, M_r$, $0<k<r$, there exist a function $x \in L_{\infty,\infty}^r(\mathbb R)$, for which these numbers are the norms of the function itself, its $k$-th and its $r$-th derivatives respectively if and only if 
$$
M_k\leq
\frac{\left\| \varphi _{r-k}\right\| }{%
\left\| \varphi _r\right\| ^{1-k/r}}%
M_0^{1-k/r}M_r^{k/r},
$$
where $\varphi _r$  ---  $r$-th periodic integral with zero mean at the period of the function $\varphi_0\left( t\right) ={\rm sgn}\sin t$. 

Another results for the case when the domain of definition of the functions is the whole real line can be found in the articles of~Rodov~[\ref{Rodov1946}, \ref{Rodov1954}],  Dzyadyk and Dubovik~[\ref{Dzyadyk1974}, \ref{Dzyadyk1975}], Babenko and Kovalenko~[\ref{BK}].

For the case when the domain of definition of the functions is semi-line~$\mathbb R_-$, the solution of Kolmogorov problem is known in the following cases:
\begin{enumerate}
\item $X=L^{r}_{\infty,\infty}({\mathbb R}_-)$, $k_1=0<k_2< k_3=r$ (partial cases follow from the results of Landau 1913~[\ref{Lan}], Matorin
1955~[\ref{Mat}]; general result follows from the work of Schoenberg and Cavaretta 1970~[\ref{ShC}]).

\item $X=L^{r}_{\infty,\infty}({\mathbb R}_-)$, $k_1=0<k_2<k_3=r-1,\;k_4=r$ (Babenko and Britvin 2002~[\ref{BB}] ).
\end{enumerate}

Let $r,m\in \mathbb{Z}_+, \;m\le r$ be given.  Denote by $L_{\infty,\infty}^{r,m}(\mathbb R_-)$ the class of functions $x \in L^{r}_{\infty,\infty}(\mathbb R_-)$, that are non-negative together with their derivatives up to the order $m$ inclusively (the derivative of the order $m$ must be non-negative almost everywhere in the case when $m=r$). We will call such functions $m$-monotone functions.

In 1951 Olovyanishnikov~[\ref{ol}] achieved the solution of Kolmogorov problem in the case $d=3$ and $X = L^{r,r-1}_{\infty,\infty}(\mathbb R_-)$.
He showed that for arbitrary $k,r\in \mathbb N $, $k<r$, and positive numbers $M_0, M_{k}, M_r$, there exists a function $x\in L^{r,r-1}_{\infty,\infty}(\mathbb R_-)$ such that 
$$
\|x\|=M_0,\;\;\|x^{(k)}\|=M_{k},\;\;\|x^{(r)}\|=M_{r},
$$
if and only if these numbers satisfy Kolmogorov type inequality
\begin{equation}\label{olov}
M_0\geq \frac{(r-k)!^{r/(r-k)}}{r!}M_k^{\frac{r}{r-k}}M_r^{-\frac{k}{r-k}}.
\end{equation}

In~[\ref{subbchern}] and independently in~[\ref{usEJA}] a generalization of this result for the case of $(r-2)$-monotone functions was achieved (moreover, in the later work case of arbitrary norms was considered).

Other results for the class of multiply monotone functions are known in the following cases:
\begin{enumerate}
\item $X=L^{r,r-2}_{\infty,\infty}({\mathbb R}_-)$ and $k_1=0< k_2< k_3=r-1,\; k_4=r$ (Yattselev 1999~[\ref{Yat}]). 
\item $X=L^{r,r}_{\infty,\infty}({\mathbb R}_-)$ and $k_1=0< k_2< k_3<k_4=r$ (V.~Babenko, Yu.~Babenko 2007~[\ref{us}]). 
\end{enumerate}

The aim of this article is to solve Kolmogorov problem for the system of positive numbers $ M_{k_1}, M_{k_2}, M_{k_3}, M_{k_4}$,  $0=
k_1<k_2<k_3 < k_4 = r$ and the class $X=L^{r,r-1}_{\infty,\infty}({\mathbb R}_-)$. Note, that the case when $k_3 = r-1$, is contained, in fact, in the work of Yattselev~[\ref{Yat}]. Hence the case $0= k_1<k_2<k_3\leq r-2$,  $k_4 = r$ is of interest.

{\bf 2. Extremal functions and some of their properties.} 

For numbers $l > 0$ and $a>b\geq 0$ let $$\varphi_1(a,b,l;t):= l\left((t+a)_+ - 2(t+b)_+\right)_+,\, t\in\mathbb R_-.$$ For $r\in \mathbb N$, $r\geq 2$ let $\varphi_r(a,b,l;t):=\int\limits_{-\infty}^t \varphi_{r-1}(a,b,l;s)ds$. Note, that  $\varphi_r(a,b,l;t)\in L^{r,r-1}_{\infty,\infty}({\mathbb R}_-)$ and $\varphi_r(a,b,l;t)=0$, $t\leq -a$. 

The following lemma holds.

\textbf{Lemma 1.} \emph{Let the function $x\in L^{r,r-1}_{\infty,\infty}({\mathbb R}_-)$ and integer numbers $0 = k_1 < k_2 <  k_{3} \leq r-2$, $k_4=r$ be given. Let numbers $l > 0$ and $a>b\geq 0$ be such, that 
\begin{equation}\label{l.extremalProperty1}
\|x^{(k_i)}\|= \|\varphi_r^{(k_i)}(a,b,l)\|,\, i=2,3,4.
\end{equation}
 Then $\|x\|\geq \|\varphi_r(a,b,l)\|$.}
 
\textbf{\emph{Proof.}} Assume the contrary, let the opposite inequality $\|x\| < \|\varphi_r(a,b,l)\|$ hold. Let $\delta(t):=x(t) - \varphi_r(a,b,l;t)$. In the virtue of the assumption and inclusion $x,\varphi_r(a,b,l;t)\in L^{r,r-1}_{\infty,\infty}({\mathbb R}_-)$  we get that $\delta(0) = x(0) - \varphi_r(a,b,l;0) = \|x\| - \|\varphi_r(a,b,l)\| < 0$. Moreover in the virtue of the definition of the functions $\varphi_r(a,b,l;t)$  $\delta(-a) \geq 0$. This means that there exists a point $-a < t_1^1<0$  such that $\delta'(t_1^1) < 0$. Moreover, $\delta'(-a) \geq 0$. This means that there exists a point $-a < t_2^1<0$ such that $\delta''(t_2^1) < 0$. Continuing with the same arguments we get that there exists a point $-a <t_{k_2}^1<0$ such that $\delta^{(k_2)}(t_{k_2}^1) < 0$. Moreover, $\delta^{(k_2)}(-a) \geq 0$, and in virtue of~\eqref{l.extremalProperty1} $\delta^{(k_2)}(0) = 0$. This means that there exist points $-a< t_{k_2+1}^1<t_{k_2+1}^2 <0$  such that $\delta^{(k_2+1)}(t_{k_2+1}^1) < 0$ and $\delta^{(k_2+1)}(t_{k_2+1}^2) > 0$. Moreover, $\delta^{(k_2+1)}(-a)  \geq 0$. Continuing with the same arguments we get that there exist points $-a< t_{k_3}^1<t_{k_3}^2 <0$  such that $\delta^{(k_3)}(t_{k_3}^1) < 0$ and $\delta^{(k_3)}(t_{k_3}^2) > 0$. Moreover,  $\delta^{(k_3)}(-a) \geq 0$ and in virtue of~\eqref{l.extremalProperty1} $\delta^{(k_3)}(0) = 0$. This means that there exist points $-a< t_{k_3+1}^1<t_{k_3+1}^2 <t_{k_3+1}^3  <0$  such that $\delta^{(k_3+1)}(t_{k_3+1}^1) < 0$, $\delta^{(k_3+1)}(t_{k_3+1}^2) > 0$ and $\delta^{(k_3+1)}(t_{k_3+1}^3) < 0$. And so on, there exist points $-a< t_{r-1}^1<t_{r-1}^2 <t_{r-1}^3  <0$  such that $\delta^{(r-1)}(t_{r-1}^1) < 0$, $\delta^{(r-1)}(t_{r-1}^2) > 0$ and $\delta^{(r-1)}(t_{r-1}^3) < 0$. But in the virtue of~\eqref{l.extremalProperty1} (with $i=4$) and the definition of the functions $\varphi_r(a,b,l;t)$ it is impossible, because at each of the intervals $(-a, -b)$ and $(-b,0)$ the function $\delta^{(r-1)}$ can have not more than one change of sign; at that from ''plus'' to ''minus'' at the first interval and from ''minus'' to ''plus'' on the second interval. Contradiction. Lemma is proved.

The following lemma can be proved analogically to lemma~1.

\textbf{Lemma 2.} \emph{Let a function $x\in L^{r,r-1}_{\infty,\infty}({\mathbb R}_-)$ and integers $0 \leq k_1 < k_2 <  k_{3} =r$ be given. Let the numbers $l > 0$ and $a> 0$ are such that
\begin{equation}\label{l.3_numbers1}
\|x^{(k_i)}\|= \|\varphi_r^{(k_i)}(a,0,l)\|,\, i=2,3.
\end{equation}
 Then $\|x^{(k_1)}\|\geq \|\varphi_r^{(k_1)}(a,0,l)\|$.}

Note, that for any function $x\in L^{r,r-1}_{\infty,\infty}({\mathbb R}_-)$ and numbers $0 \leq k_1 < k_2 <  k_{3} =r$ parameters  $a,l> 0$ can be chosen in such a way, that the equalities~\eqref{l.3_numbers1} hold. Moreover, since $\varphi_r(a,0,l)=\frac{l}{r!}(t+a)^r_+$, we have that
\begin{equation}\label{1}
\left\|\varphi_r^{(k_1)}(a,0,l)\right\|= \frac{(r-k_2)!^{\frac{r-k_1}{r-k_2}}}{(r-k_1)!}\left\|\varphi_r^{(k_2)}(a,0,l)\right\|^{\frac{r-k_1}{r-k_2}}\left\|\varphi_r^{(r)}(a,0,l)\right\|^{\frac{k_1-k_2}{r-k_2}}.
\end{equation}

Hence the following lemma holds.

\textbf{Lemma 3.} \emph{Let the function $x\in L^{r,r-1}_{\infty,\infty}({\mathbb R}_-)$ and integers $0 \leq k_1 < k_2 <  k_{3} =r$ be given. 
 Then $$\left\|x^{(k_1)}\right\|\geq \frac{(r-k_2)!^{\frac{r-k_1}{r-k_2}}}{(r-k_1)!}\left\|x^{(k_2)}\right\|^{\frac{r-k_1}{r-k_2}}\left\|x^{(r)}\right\|^{\frac{k_1-k_2}{r-k_2}}.$$}

Lemma~3 is generalization of Olovyanishnikov inequality~\eqref{olov} and, in fact, is contained in~[\ref{usEJA}].

\textbf{Lemma 4.} \emph{Let integers $0 \leq k_1 < k_2\leq r-2$, $ k_{3} =r$ and positive numbers $M_{k_1},M_{k_2},M_{r}$ such that the following inequality  holds
$$M_{k_1} \geq \frac{(r-k_2)!^{\frac{r-k_1}{r-k_2}}}{(r-k_1)!}M_{k_2}^{\frac{r-k_1}{r-k_2}}M_r^{\frac{k_1-k_2}{r-k_2}}$$ 
be given.
 Then there exist numbers $l>0$, $a>b\geq 0$, such that the following equalities hold:
\begin{equation}\label{l.splineExistance1}
\|\varphi_r^{(k_i)}(a,b,l)\|=M_{k_i},\, i=1,2,3.
\end{equation}}

\textbf{\emph{Proof.}}
Note, that in virtue of the definition of  $\varphi_r(a,b,l)$ we have $\|\varphi_r^{(r)}(a,b,l)\| = l$. Hence below we can count that $M_r=l=1$ and instead of $\varphi_r(a,b,1)$ we will write $\varphi_r(a,b)$. 

For each $b\geq 0$ there exists a number $a=a(b)$ such that 
\begin{equation}\label{l.splineExistance2}
\|\varphi_r^{(k_2)}(a(b),b)\| = M_{k_2}.
\end{equation}
 Really, for fixed $b$, $\psi(a):=\|\varphi_r^{(k_2)}(a,b)\| $ is a continuous function of variable $a$. Moreover, $\psi(b)=0$ and $\psi(a)\to\infty$ when $a\to\infty$. This means that there exists a number $a= a(b)$ such that $\psi(a(b))=M_{k_2}$. Hence on the interval $[0,\infty)$ we defined a function $a(b)$ such that for all $b\geq 0$ the equality~\eqref{l.splineExistance2} holds. At the same time the function $a(b)$ is continuous.
 
 We will show, that there exists a number $b\geq 0$, such that 
 \begin{equation}\label{l.splineExistance3}
\|\varphi_r^{(k_1)}(a(b),b)\| = M_{k_1}.
\end{equation}

Let $\eta(b):=\|\varphi_r^{(k_1)}(a(b),b)\|$. By the conditions of the lemma $\eta(0)\leq M_{k_1}$. We will now show that 
\begin{equation}\label{l.splineExistance3.1}
\eta(b)\to\infty,\, b\to\infty.
\end{equation}
Since the equality~\eqref{l.splineExistance2} holds, in virtue of the definition of the functions $\varphi(a,b,l;t)$ the value $a(b)-b$ is bounded. Hence
\begin{equation}\label{l.splineExistance4}
2b-a(b)\to\infty,\, b\to\infty.
\end{equation}
From the definition of the functions $\varphi(a,b,l;t)$  it follows, that the restriction $p(b;t)$ of the function $\varphi(a(b),b;t)$ to the segment $[a(b)-2b,0]$ is polynomial of degree  $r-2$. Moreover, 
\begin{equation}\label{l.splineExistance5}
 \max\limits_{t\in [a(b)-2b,0]}|p^{(k_2)}(b;t)| = p^{(k_2)}(b;0) = M_{k_2}>0. 
\end{equation} 

Applying Markov inequality for algebraic polynomials (see.~[\ref{MarkovVA}],~[\ref{MarkovAA}], see also chapter~4 in~[\ref{KBL}]), since~\eqref{l.splineExistance4} and~\eqref{l.splineExistance5} we get that 
$$\max\limits_{t\in [a(b)-2b,0]}|p^{(k_1)}(b;t)|\to\infty,\,b\to\infty.$$
Hence we get the truth of~\eqref{l.splineExistance3.1}. Thus there exist numbers $l>0$, $a>b \geq 0$ such that the equalities~\eqref{l.splineExistance1} hold. Lemma is proved.

{\bf Remark. } {\it Under the conditions of lemma~4 let $$\Phi(M_{k_1}, M_{k_2}, M_{k_3};t) := \varphi(a,b,l;t),$$ where the numbers $a,b,l$ are chosen in such a way that the equalities~\eqref{l.splineExistance1} hold.}

{\bf 3. Solution of Kolmogorov problem.}

\textbf{Theorem 1.} \emph{Let integers $0 =k_1\leq k_2 < k_3\leq r-2$, $ k_4 =r$ and positive numbers $M_{k_1},M_{k_2},M_{k_3}, M_{k_4}$ be given. There exist a function $x\in L^{r,r-1}_{\infty,\infty}({\mathbb R}_-)$ such that $$\|x^{(k_i)}\|= M_{k_i},\,i=1,2,3,4$$ if and only if the following inequalities hold:
 $$M_{k_2} \geq \frac{(r-k_3)!^{(r-k_2)/(r-k_3)}}{(r-k_2)!}M_{k_3}^{\frac{r-k_2}{r-k_3}}M_r^{\frac{k_2-k_3}{r-k_3}}$$
 and
 $$M_{k_1} \geq \|\Phi(M_{k_2}, M_{k_3}, M_{r})\|.$$}
 \textbf{\emph{Proof.}}
 The necessity of inequalities above follow from lemma~3 and lemma~1. If the inequalities above hold then the function $x(t) := \Phi(M_{k_2}, M_{k_3}, M_{r};t) + M_{k_1} - \|\Phi(M_{k_2}, M_{k_3}, M_{r})\|$ is desired. Theorem is proved.

\vskip 3.5mm
\footnotesize
\begin{enumerate}

\item\label{Kolmogorov}
{\it Kolmogorov, A. N.}: { Selected works of A. N. Kolmogorov. Vol. I. Mathematics and mechanics.} Translation: Mathematics and its Applications (Soviet Series), 25. Kluwer Academic Publishers Group, Dordrecht. --- 1991. p 252--263.

\item\label{had}
{\it Hadamard J.} {Sur le maximum d'une fonction et de ses derivees} // C. R. Soc. Math. France. --- 1914. --- {\bf 41}. --- P. 68--72.

\item\label{shilov}
{\it Shilov G. E.} {On inequalities between derivatives} // In the book '' Sbornik rabot studencheskih nauchnyh kruzhkov Mosc. Univ.'' --- 1937. --- {\bf 1}. --- p. 17--27 (in Russian).

\item\label{Rodov1946}
{\it Rodov A. M. } {Dependence between upper bounds of arbitrary functions of real variable} // Izv. AN USSR. Ser. Math. ---
1946. --- {\bf 10}. --- C. 257--270 (in Russian).

\item\label{Rodov1954}
{\it Rodov A. M. }  {Sufficient conditions of the existence of a function of real variable with prescribed upper bounds of moduli of the function itself and its five consecutive derivatives} // Uchenye Zapiski Belorus. Univ. --- 1954. --- {\bf 19}. --- p. 65--72 (in Russian).

\item\label{Dzyadyk1974}
{\it  Dzyadyk V. K., Dubovik V. A.} {On problem of A.~N.~Kolmogorov about dependence between upper bounds of the derivatives of real value functions given on the whole line} // Ukr. Math. Journ. --- 1974. --- {\bf 26} (3). --- p. 300--317 (in Russian).

\item\label{Dzyadyk1975}
{\it Dzyadyk V. K., Dubovik V. A.} {On inequalities of A.~N.~Kolmogorov about dependence between upper bounds of the
derivatives of real value functions given on the whole line} // Ukr. Math. Journ. --- 1975. --- {\bf 27} (3). --- p. 291--299 (in Russian).

\item\label{BK}
{\it Babenko V. F., Kovalenko O. V.} On the dependence of the norm of a function on the norms of its derivatives of orders $k$, $r - 2$, and $r$, $0 < k < r - 2$ // Ukr. Math. Journ. --- 2012. --- {\bf 64} (5). --- p. 597--603.

\item\label{Lan}
{\it Landau E. } { Einige Ungleichungen fur zweimal differenzierbare Funktion} // Proc. London Math. Soc. --- 1913. --- {\bf 13}. --- P. 43--49.

\item\label{Mat}
{\it Matorin, A. P. } { On inequalities between the maxima of the absolute values of a function and its derivatives on a half-line}, Ukr. Math. Journ.,  {\bf 7}. --- 1951. --- p. 262--266 (in Russian).

\item\label{ShC}
{\it Schoenberg I. J., Cavaretta A.}  Solution of Landau's problem, concerning higher derivatives on half line // Proc. of Conference on Approximation theory. Varna 1970. --- P. 297--308.

\item\label{BB}
{\it Babenko, V. F., Britvin, Y. E.}  On Kolmogorov's problem about existence of a function with given norms of its derivatives // East J. Approx. --- 2002. --- {\bf 8}, №1. --- {P. 95--100.}

\item\label{ol}
{\it Olovyanishnikov V. M.} To the question on inequalities between upper bounds of consecutive derivatives on a half-line // Uspehi mat. nauk --- 1951. --- {\bf 6}(2)(42). --- {С. 167–170}

\item\label{subbchern}
{\it Subbotin, Yu. N., Chernyh, N. I.} Inequalities for derivatives of monotone functions // Approximation of functions. Theoretical and applied aspects. Coll. papers.--- 2003. --- {p. 199-211} (in Russian).

\item\label{usEJA}
{\it Babenko V., Babenko Yu.} The Kolmogorov Inequalities for Multiply Monotone Functions Defined on a Half-line // East J. Approx. --- 2005. ---{\bf 11}, №2. --- {P. 169--186.}

\item\label{Yat}
{\it Yattselev  M. L. }  {Inequality between four upper bounds of consecutive derivatives on a half-line } // Visn. DGU.
Mathematics, {\bf 4}. --- 1998. --- p. 106--111 (in Russian).

\item\label{us}
{\it Babenko V., Babenko Yu.}   On the Kolmogorov's problem for the upper bounds of four consecutive derivatives of a multiply monotone function // Constr. Approx. --- 2007. --- {\bf 26}, №1. --- P. 83--92.

\item\label{MarkovVA} {\it Markov V. A.} On functions of least deviation from zero in a given interval // SPb --- 1892.  St. Petersburg (in Russian).

\item\label{MarkovAA} {\it Markov A. A.} On a question by D. I. Mendeleev. Zap. Imp. Akad. Nauk SPb. --- 1989. --- {\bf 62}. --- p. 1--24 (in Russian).

\item\label{KBL} {\it Korneichuk N. P., Ligun A. A., Babenko V. F. } Extremal properties of polynomials and splines // Kyiv: Naukova Dumka. --- 1992. --- 304 p.

\end{enumerate}

\label{end}
\end{document}


\begin{center}
{\bf \Large Задача Колмогорова на классе кратно монотонных функций}

\bigskip

О. В. Коваленко

\bigskip

Днепропетровский национальный университет им. О. Гончара

\end{center}
\vskip 2mm
 {\small \it \noindent Найдены необходимые и достаточные условия на положительные числа
$M_{k_1}, M_{k_2},  M_{k_3}, M_{k_4}$, $0 = k_1 < k_2<k_2\leq r-2$, $k_4=r$, для того, чтобы гарантировать существование $r-1$-кратно монотонной функции определенной на отрицательной полуоси и такой, что $\|x^{(k_i)}\| = M_{k_i}$, $i=1,2,3,4$. \hfill}

 \normalsize \vskip 3mm

{\bf 1. Обозначения. Постановка задачи. Известные результаты.} 
Пусть область $G$ обозначает действительную ось  $\mathbb R=(-\infty , \infty )$ или неотрицательную полуось $\mathbb R_-=(-\infty , 0]$. Пусть  $L_{\infty}(G)$ обозначает пространство измеримых существенно ограниченных функций  $x:G\to \mathbb R$ с обычной нормой $\| \cdot \|=\|
\cdot \|_{L_{\infty}(G)}$. Для $r\in \mathbb N$ обозначим через $L^r_{\infty} (G)$  пространство функций $x:G\to \mathbb R$, которые имеют локально абсолютно непрерывную производную порядка $r-1$, $x^{(0)}=x$, и таких, что $x ^ {(r)} \in L_{\infty} (G)$. Положим $L^r_{\infty,\infty}(G)=L^r_{\infty}(G)\cap L_{\infty}(G)$. 

Для действительного $t\in \mathbb R$ положим $t_+ := \max\{t,0\}$.

А.Н. Колмогоров (см.~{[\ref{Kolmogorov}]}) сформулировал следущую задачу:

{\bf Задача Колмогорова}

{\it Пусть задан некоторый класс функций $X\subset L^r_{\infty,\infty}(G)$ и произвольная система $d$ целых чисел
$0 = k_1 < k_2 < ... < k_{d}=r.$
Найти необходимые и достаточные условия на систему положительных чисел 
$$M_{k_1}, M_{k_2}, ..., M_{k_{d}}$$
для того, чтобы гарантировать существование функции $x \in X$
 такой, что $$\|x^{(k_i)}\|=M_{k_i},\; i=1,...,d.$$

}

 В~{[\ref{Kolmogorov}] А.Н. Колмогоров решил эту задачу в случае $d=3$, $X= L^r_{\infty,\infty}(\mathbb R)$ (частные случаи следуют из работ Адамара~[\ref{had}] и Шилова~[\ref{shilov}]). 
Он показал, что для трех положительных чисел $M_0, M_k, M_r$, $0<k<r$, существует функция $x \in L_{\infty,\infty}^r(\mathbb R)$, для которой эти числа являются нормами функии, ее $k$-й и ее $r$-й производных соответственно, тогда и только тогда, когда 
$$
M_k\leq
\frac{\left\| \varphi _{r-k}\right\| }{%
\left\| \varphi _r\right\| ^{1-k/r}}%
M_0^{1-k/r}M_r^{k/r},
$$
где $\varphi _r$  ---  $r$-я периодическая первообразная со средним значением ноль на периоде функции $\varphi_0\left( t\right) ={\rm sgn}\sin t$. 

Другие результаты в случае, когда областью определения фукнций является вся дейст\-вительная ось содержатся в статьях Родова~[\ref{Rodov1946}, \ref{Rodov1954}],  Дзядыка и Дубовика~[\ref{Dzyadyk1974}, \ref{Dzyadyk1975}], Бабенко и Коваленко~[\ref{BK}].

Остановимся подробнее на случае, когда областью определения является полуось $\mathbb R_-$. Решение задачи Колмогорова известно в следующих случаях
\begin{enumerate}
\item $X=L^{r}_{\infty,\infty}({\mathbb R}_-)$, $k_1=0<k_2< k_3=r$ (частные случаи следуют из результатов  Ландау 1913~[\ref{Lan}], Маторина
1955~[\ref{Mat}]; общий результат следует из работы Шенберга и Каваретта 1970,~[\ref{ShC}]).

\item $X=L^{r}_{\infty,\infty}({\mathbb R}_-)$, $k_1=0<k_2<k_3=r-1,\;k_4=r$ (Бабенко и  Бритвин 2002~[\ref{BB}] ).
\end{enumerate}

Пусть заданы $r,m\in \mathbb{Z}_+, \;m\le r$. Обозначим через $L_{\infty,\infty}^{r,m}(\mathbb R_-)$ класс функций  $x \in L^{r}_{\infty,\infty}(\mathbb R_-)$, которые неотрицательны вместе со всеми своими производными до порядка $m$ включительно (производная порядка $m$ должна быть неотрицательной почти всюду в случае $m=r$). Мы будем называть этот класс классом $m$-кратно монотонных функций.

В 1951 году Оловянишников~[\ref{ol}] получил решение задачи Колмогорова в случае $d=3$ и $X = L^{r,r-1}_{\infty,\infty}(\mathbb R_-)$.
Он показал, что для произвольных $k,r\in \mathbb N $, $k<r$, и положительных чисел $M_0, M_{k}, M_r$, существует функция $x\in L^{r,r-1}_{\infty,\infty}(\mathbb R_-)$ такая, что 
$$
\|x\|=M_0,\;\;\|x^{(k)}\|=M_{k},\;\;\|x^{(r)}\|=M_{r},
$$
тогда и только тогда, когда эти три числа удовлетворяют неравенству колмогоровского типа
\begin{equation}\label{olov}
M_0\geq \frac{(r-k)!^{r/(r-k)}}{r!}M_k^{\frac{r}{r-k}}M_r^{-\frac{k}{r-k}}.
\end{equation}

В~[\ref{subbchern}] и независимо в~[\ref{usEJA}] было получено обобщение этого результата на случай $(r-2)$-кратно монотонных функций (более того, в последней работе рассматривался случай произвольных норм).

Другие результаты на классах кратно монотонных функций известны в следующих случаях:
\begin{enumerate}
\item $X=L^{r,r-2}_{\infty,\infty}({\mathbb R}_-)$ и $k_1=0< k_2< k_3=r-1,\; k_4=r$ (Ятцелев 1999~[\ref{Yat}]). 
\item $X=L^{r,r}_{\infty,\infty}({\mathbb R}_-)$ и $k_1=0< k_2< k_3<k_4=r$ (В.~Бабенко, Ю.~Бабенко 2007~[\ref{us}]). 
\end{enumerate}

Цель данной статьи -- решение задачи Колмогорова для системы
положи\-тельных чисел $ M_{k_1}, M_{k_2}, M_{k_3}, M_{k_4}$,  $0=
k_1<k_2<k_3 < k_4 = r$ и класса $X=L^{r,r-1}_{\infty,\infty}({\mathbb R}_-)$. Отметим, что случай, когда $k_3 = r-1$, по сути содержится в работе Ятцелева~[\ref{Yat}]. Поэтому интерес представляет рассмотрение случая $0= k_1<k_2<k_3\leq r-2$,  $k_4 = r$.

{\bf 2. Экстремальные фунции и некоторые их свойства.} 

Для чисел $l > 0$ и $a>b\geq 0$ положим $$\varphi_1(a,b,l;t):= l\left((t+a)_+ - 2(t+b)_+\right)_+,\, t\in\mathbb R_-.$$ Для $r\in \mathbb N$, $r\geq 2$ положим $\varphi_r(a,b,l;t):=\int\limits_{-\infty}^t \varphi_{r-1}(a,b,l;s)ds$. Отметим, что  $\varphi_r(a,b,l;t)\in L^{r,r-1}_{\infty,\infty}({\mathbb R}_-)$ и $\varphi_r(a,b,l;t)=0$, $t\leq -a$. 

Справедлива следующая лемма.

\textbf{Лемма 1.} \emph{Пусть заданы функция $x\in L^{r,r-1}_{\infty,\infty}({\mathbb R}_-)$ и целые числа $0 = k_1 < k_2 <  k_{3} \leq r-2$, $k_4=r$. Пусть числа $l > 0$ и $a>b\geq 0$ таковы, что 
\begin{equation}\label{l.extremalProperty1}
\|x^{(k_i)}\|= \|\varphi_r^{(k_i)}(a,b,l)\|,\, i=2,3,4.
\end{equation}
 Тогда $\|x\|\geq \|\varphi_r(a,b,l)\|$.}
 
\textbf{\emph{Доказательство.}} Предположим противное, пусть имеет место неравенство $\|x\| < \|\varphi_r(a,b,l)\|$. Положим $\delta(t):=x(t) - \varphi_r(a,b,l;t)$. В силу предположения и принадлежности $x,\varphi_r(a,b,l;t)\in L^{r,r-1}_{\infty,\infty}({\mathbb R}_-)$  мы получаем, что $\delta(0) = x(0) - \varphi_r(a,b,l;0) = \|x\| - \|\varphi_r(a,b,l)\| < 0$. Кроме того, в силу определения фукнций $\varphi_r(a,b,l;t)$  $\delta(-a) \geq 0$. Это значит, что существует точка $-a < t_1^1<0$  такая, что $\delta'(t_1^1) < 0$. Кроме того, $\delta'(-a) \geq 0$. Это значит, что существует точка $-a < t_2^1<0$ такая, что $\delta''(t_2^1) < 0$. Повторяя аналогичные рассуждения получим, что существует точка $-a <t_{k_2}^1<0$ такая, что $\delta^{(k_2)}(t_{k_2}^1) < 0$. Кроме того, $\delta^{(k_2)}(-a) \geq 0$ и в силу~\eqref{l.extremalProperty1} $\delta^{(k_2)}(0) = 0$. Это значит, что существуют точки $-a< t_{k_2+1}^1<t_{k_2+1}^2 <0$  такие, что $\delta^{(k_2+1)}(t_{k_2+1}^1) < 0$ и $\delta^{(k_2+1)}(t_{k_2+1}^2) > 0$. Кроме того, $\delta^{(k_2+1)}(-a)  \geq 0$. Повторяя аналогичные рассуждения мы получим, что существуют точки $-a< t_{k_3}^1<t_{k_3}^2 <0$  такие, что $\delta^{(k_3)}(t_{k_3}^1) < 0$ и $\delta^{(k_3)}(t_{k_3}^2) > 0$. Кроме того,  $\delta^{(k_3)}(-a) \geq 0$ и в силу~\eqref{l.extremalProperty1} $\delta^{(k_3)}(0) = 0$. Это значит, что существуют точки $-a< t_{k_3+1}^1<t_{k_3+1}^2 <t_{k_3+1}^3  <0$  такие, что $\delta^{(k_3+1)}(t_{k_3+1}^1) < 0$, $\delta^{(k_3+1)}(t_{k_3+1}^2) > 0$ и $\delta^{(k_3+1)}(t_{k_3+1}^3) < 0$. И так далее, существуют точки $-a< t_{r-1}^1<t_{r-1}^2 <t_{r-1}^3  <0$  такие, что $\delta^{(r-1)}(t_{r-1}^1) < 0$, $\delta^{(r-1)}(t_{r-1}^2) > 0$ и $\delta^{(r-1)}(t_{r-1}^3) < 0$.  Однако в силу~\eqref{l.extremalProperty1} (при $i=4$) и определения функций $\varphi_r(a,b,l;t)$ это невозможно, поскольку на каждом из интервалов $(-a, -b)$ и $(-b,0)$ функция $\delta^{(r-1)}$ может иметь не более одной перемены знака, причем с ''плюс'' на ''минус'' на первом и с ''минус'' на ''плюс'' на втором. Пришли к противоречию. Лемма доказана.

Следущая лемма доказывается аналогично лемме~1.

\textbf{Лемма 2.} \emph{Пусть заданы функция $x\in L^{r,r-1}_{\infty,\infty}({\mathbb R}_-)$ и целые числа $0 \leq k_1 < k_2 <  k_{3} =r$. Пусть числа $l > 0$ и $a> 0$ таковы, что 
\begin{equation}\label{l.3_numbers1}
\|x^{(k_i)}\|= \|\varphi_r^{(k_i)}(a,0,l)\|,\, i=2,3.
\end{equation}
 Тогда $\|x^{(k_1)}\|\geq \|\varphi_r^{(k_1)}(a,0,l)\|$.}

Отметим, что для любой функции $x\in L^{r,r-1}_{\infty,\infty}({\mathbb R}_-)$ и целых чисел $0 \leq k_1 < k_2 <  k_{3} =r$ можно выбрать параметры $a,l> 0$ так, чтобы выполнялись равенства~\eqref{l.3_numbers1}. Кроме того, так как $\varphi_r(a,0,l)=\frac{l}{r!}(t+a)^r_+$, то  справедливо равенство 
\begin{equation}\label{1}
\left\|\varphi_r^{(k_1)}(a,0,l)\right\|= \frac{(r-k_2)!^{\frac{r-k_1}{r-k_2}}}{(r-k_1)!}\left\|\varphi_r^{(k_2)}(a,0,l)\right\|^{\frac{r-k_1}{r-k_2}}\left\|\varphi_r^{(r)}(a,0,l)\right\|^{\frac{k_1-k_2}{r-k_2}}.
\end{equation}

Поэтому справедлива следующая лемма.

\textbf{Лемма 3.} \emph{Пусть заданы функция $x\in L^{r,r-1}_{\infty,\infty}({\mathbb R}_-)$ и целые числа $0 \leq k_1 < k_2 <  k_{3} =r$. 
 Тогда $$\left\|x^{(k_1)}\right\|\geq \frac{(r-k_2)!^{\frac{r-k_1}{r-k_2}}}{(r-k_1)!}\left\|x^{(k_2)}\right\|^{\frac{r-k_1}{r-k_2}}\left\|x^{(r)}\right\|^{\frac{k_1-k_2}{r-k_2}}.$$}

Лемма~3 является обобщением неравенства Оловянишникова~\eqref{olov} и, по сути, содержится в~[\ref{usEJA}].

\textbf{Лемма 4.} \emph{Пусть заданы целые числа $0 \leq k_1 < k_2\leq r-2$, $ k_{3} =r$ и положительные числа $M_{k_1},M_{k_2},M_{r}$ такие, что выполняется неравенство 
$$M_{k_1} \geq \frac{(r-k_2)!^{\frac{r-k_1}{r-k_2}}}{(r-k_1)!}M_{k_2}^{\frac{r-k_1}{r-k_2}}M_r^{\frac{k_1-k_2}{r-k_2}}.$$
 Тогда существуют такие числа $l>0$, $a>b\geq 0$, что выполняются равенства 
\begin{equation}\label{l.splineExistance1}
\|\varphi_r^{(k_i)}(a,b,l)\|=M_{k_i},\, i=1,2,3.
\end{equation}}

\textbf{\emph{Доказательство.}}
Отметим, что в силу определения функций $\varphi_r(a,b,l)$ имеем $\|\varphi_r^{(r)}(a,b,l)\| = l$. Поэтому далее можем считать, что $M_r=l=1$ и вместо $\varphi_r(a,b,1)$ будем писать $\varphi_r(a,b)$. 

Для каждого $b\geq 0$ существует число $a=a(b)$ такое, что 
\begin{equation}\label{l.splineExistance2}
\|\varphi_r^{(k_2)}(a(b),b)\| = M_{k_2}.
\end{equation}
 Действительно, при фиксированном $b$  $\psi(a):=\|\varphi_r^{(k_2)}(a,b)\| $ есть непрерывная функция переменного $a$. Кроме того, $\psi(b)=0$ и $\psi(a)\to\infty$ при $a\to\infty$. Это значит, что существует число $a= a(b)$ такое, что $\psi(a(b))=M_{k_2}$. Таким образом на промежутке $[0,\infty)$ мы определили функцию $a(b)$ такую, что для всех $b\geq 0$ выполняется равенство~\eqref{l.splineExistance2}. При этом фукнция $a(b)$ является непрерывной.
 
 Покажем, что существует $b\geq 0$, такое, что выполняется 
 \begin{equation}\label{l.splineExistance3}
\|\varphi_r^{(k_1)}(a(b),b)\| = M_{k_1}.
\end{equation}

Положим $\eta(b):=\|\varphi_r^{(k_1)}(a(b),b)\|$. По условию леммы $\eta(0)\leq M_{k_1}$. Покажем, что 
\begin{equation}\label{l.splineExistance3.1}
\eta(b)\to\infty,\, b\to\infty.
\end{equation}
Поскольку выполняется равенство~\eqref{l.splineExistance2}, то в силу определения функций $\varphi(a,b,l;t)$ величина $a(b)-b$ ограничена. Поэтому 
\begin{equation}\label{l.splineExistance4}
2b-a(b)\to\infty,\, b\to\infty.
\end{equation}
 Из определения функций $\varphi(a,b,l;t)$ следует, что сужение $p(b;t)$ функции $\varphi(a(b),b;t)$ на отрезок $[a(b)-2b,0]$ есть полином степени  $r-2$. Кроме того, 
\begin{equation}\label{l.splineExistance5}
 \max\limits_{t\in [a(b)-2b,0]}|p^{(k_2)}(b;t)| = p^{(k_2)}(b;0) = M_{k_2}>0. 
\end{equation} 

Теперь применяя неравенство Маркова для алгебраических полиномов (см.~[\ref{MarkovVA}],~[\ref{MarkovAA}], см. также гл.~4 в монографии~[\ref{KBL}]) и учитывая~\eqref{l.splineExistance4} и~\eqref{l.splineExistance5} получаем, что 
$$\max\limits_{t\in [a(b)-2b,0]}|p^{(k_1)}(b;t)|\to\infty,\,b\to\infty.$$
Отсюда получаем справедливость~\eqref{l.splineExistance3.1}. Таким образом существуют числа $l>0$, $a>b \geq 0$ такие, что выполняются равенства~\eqref{l.splineExistance1}. Лемма доказана.

{\bf Замечание. } {\it При выполнении условий леммы~4 положим $$\Phi(M_{k_1}, M_{k_2}, M_{k_3};t) = \varphi(a,b,l;t),$$ где числа $a,b,l$ выбраны так, что  выполняются равенства~\eqref{l.splineExistance1}.}

{\bf 3. Решение задачи Колмогорова.}

\textbf{Теорема 1.} \emph{Пусть заданы целые числа $0 =k_1\leq k_2 < k_3\leq r-2$, $ k_4 =r$ и положительные числа $M_{k_1},M_{k_2},M_{k_3}, M_{k_4}$. Существует функция $x\in L^{r,r-1}_{\infty,\infty}({\mathbb R}_-)$ такая, что $$\|x^{(k_i)}\|= M_{k_i},\,i=1,2,3,4$$ тогда и только тогда, когда выполняются следующие неравенства
 $$M_{k_2} \geq \frac{(r-k_3)!^{(r-k_2)/(r-k_3)}}{(r-k_2)!}M_{k_3}^{\frac{r-k_2}{r-k_3}}M_r^{\frac{k_2-k_3}{r-k_3}}$$
 и
 $$M_{k_1} \geq \|\Phi(M_{k_2}, M_{k_3}, M_{r})\|.$$}
 \textbf{\emph{Доказательство.}}
 Необходимость указанных неравенств следует из лемм~3 и~1. В случае выполнения указанных неравенств функция $x(t) := \Phi(M_{k_2}, M_{k_3}, M_{r};t) + M_{k_1} - \|\Phi(M_{k_2}, M_{k_3}, M_{r})\|$ является искомой. Теорема доказана.

\vskip 3.5mm
\footnotesize
\begin{enumerate}

\item\label{Kolmogorov}
{\it Колмогоров А.Н.} {О неравенствах между верхними гранями
последовательных производных произвольной функции на бесконечном
интервале.} // В кн. А.Н. Колмогоров, Избранные труды, Математика и
механика, М. Наука, 1985. --- С. 252 -- 263.

\item\label{had}
{\it Hadamard J.} {Sur le maximum d'une fonction et de ses derivees} // C. R. Soc. Math. France. --- 1914. --- {\bf 41}. --- P. 68--72.

\item\label{shilov}
{\it Шилов Г. Е.} {О неравенствах между производными }// Сборник работ студ. науч. кружков  МГУ. --- 1937. --- {\bf 1}. --- C. 17--27.

\item\label{Rodov1946}
{\it Родов А.М. } {Зависимость между верхними гранями производных функций действительного переменного} // Изв. АН СССР. Сер. Мат.---
1946. --- {\bf 10}. --- C. 257--270.

\item\label{Rodov1954}
{\it Родов А.М.}  {Достаточные условия существования функции действительного переменного с заданными верхними гранями модулей
самой функции и её пяти последовательных производных.} // Ученые записки БГУ имени В.И.Ленина. Серия физико-математическая. --- 1954.
--- {\bf 19}. --- C. 65--72.

\item\label{Dzyadyk1974}
{\it Дзядык В. К., Дубовик В. А.} {К проблеме А.Н.~Колмогорова о зависимостях между верхними гранями производных вещественных
функций, заданных на всей оси} // Укр. мат. журн. --- 1974. --- {\bf 26}, №3. --- C. 300--317.

\item\label{Dzyadyk1975}
{\it Дзядык В. К., Дубовик В. А.} {К неравенствам А.Н.~Колмогорова о зависимостях между верхними гранями производных вещественных функций, заданных на всей оси} // Укр. мат. журн. --- 1975. --- {\bf 27}, №3. --- C. 291--299.

\item\label{BK}
{\it Бабенко В.Ф., Коваленко О.В.} О зависимости между нормой функции и нормами ее производных порядка $k$, $r-2$  и $r$, $0<k<r-2$ // Укр. мат. журн. --- 2012. --- {\bf 64}, №5. --- С. 597--603.

\item\label{Lan}
{\it Landau E. } { Einige Ungleichungen fur zweimal
differenzierbare Funktion} // Proc. London Math. Soc. --- 1913. --- {\bf 13}. --- P. 43--49.

\item\label{Mat}
{\it Маторин А.П.} О неравенствах между максимумами абсолютных значений функций и ее производных на полуоси // Укр. мат. журн. --- 1955. --- {\bf 7}. --- С. 262–266.

\item\label{ShC}
{\it Schoenberg I. J., Cavaretta A.}  Solution of Landau's problem, concerning higher derivatives on half line // Proc. of Conference on Approximation theory. Varna 1970. --- P. 297--308.

\item\label{BB}
{\it Babenko, V. F., Britvin, Y. E.}  On Kolmogorov's problem about existence of a function with given norms of its derivatives // East J. Approx. --- 2002. --- {\bf 8}, №1. --- {P. 95--100.}

\item\label{ol}
{\it Оловянишников В.М.} К вопросу о неравенствах между верхними гранями последовательных производных на полупрямой // Успехи Мат. Наук. --- 1951. --- {\bf 6}(2)(42). --- {С. 167–170}

\item\label{subbchern}
{\it Субботин Ю.Н., Черных Н.И.} Неравенства для производных монотонных функций // Приближение функций. Теорет. и прикл. аспекты: Сб. ст., посвящ. памяти проф. А.В. Ефимова.--- МИЭТ.--- 2003. --- {C. 199-211.}

\item\label{usEJA}
{\it Babenko V., Babenko Yu.} The Kolmogorov Inequalities for Multiply Monotone Functions Defined on a Half-line // East J. Approx. --- 2005. ---{\bf 11}, №2. --- {P. 169--186.}

\item\label{Yat}
{\it Ятцелев М.Л.} Неравенство между четырьмя верхними гранями последовательных производных на полупрямой // Вiсник Днiпропетровського унiверситету. Математика. --- 1998. --- {\bf 4}. --- С. 106--111.

\item\label{us}
{\it Babenko V., Babenko Yu.}   On the Kolmogorov's problem for the upper bounds of four consecutive derivatives of a multiply monotone function // Constr. Approx. --- 2007. --- {\bf 26}, №1. --- P. 83--92.

\item\label{MarkovVA} {\it Марков В.А.}  О функциях, наименее уклоняющихся от нуля в данном промежутке // СПб. --- 1892. --- Санктъ-Пб.

\item\label{MarkovAA} {\it Марков А.А.} Об одном вопросе Д.И. Менделеева // Изв. АН. --- СПб. --- 1989. --- {\bf 62}. --- C.1--24.

\item\label{KBL} {\it Корнейчук Н.П., Бабенко В.Ф., Лигун А.А.} Экстремальные свойства полиномов и сплайнов. --- Киев: Наук. думка, 1992. --- 304 с.

\end{enumerate}

\label{end}